\documentclass[reqno, 12pt]{amsart}
\usepackage[letterpaper,hmargin=1in,vmargin=1in]{geometry}
\usepackage{amsmath, amssymb, amsthm, verbatim, url}

\renewcommand{\phi}{\varphi}

\newtheorem{maintheo}{{\sc Theorem}}

\newenvironment{rem}{\medskip\noindent{\it Remark:\/} }{\medskip}

\title[Lower bounds for volumes of nodal sets]
{Lower bounds for volumes of nodal sets: an improvement of a result of Sogge-Zelditch}

\author{Hamid Hezari }
\address{Department of Mathematics, MIT, Cambridge, MA 02139, USA} \email{hezari@math.mit.edu}

\author{Zuoqin Wang }
\address{Department of Mathematics, University of Michigan, Ann Arbor, MI 48109, USA} \email{wangzuoq@umich.edu}

\thanks{The first author is partially supported by the NSF
grant DMS-0969745}

\date{\today}

\begin{document}

\begin{abstract} We use the Dong-Sogge-Zelditch formula to obtain a lower bound for the volume of the nodal sets of eigenfunctions. Our result improves the recent results of Sogge-Zelditch \cite{SZ} and in dimensions $n\leq 5$ gives a new proof for the lower bounds of Colding-Minicozzi \cite{CM}.
\end{abstract}

\maketitle

\section{Introduction}
Let $(M,g)$ be a $C^\infty$ boundaryless compact Riemannian manifold of dimension $n$ and $\Delta$ be the Laplace-Beltrami operator. The eigenfunctions $\phi_\lambda$ of $\Delta$ are non zero solutions to
$$ -\Delta \phi_\lambda = \lambda \phi_\lambda.$$ Throughout the paper we assume that the eigenfunctions are normalized so that $||\phi_\lambda||_{L^2}=1$. Since $M$ is compact the spectrum of $\Delta$ is discrete and appears as an increasing sequence
$$ \lambda: \; 0=\lambda_0<\lambda_1 \leq \lambda_2 \dots \to \infty\; ,$$ where each eigenvalue has finite multiplicity. In this paper we are interested in the nodal sets
$$Z_{\phi_\lambda}=\{ x \in M;  \; \phi_\lambda(x)=0\}. $$ In particular we would like to find lower bounds for the $(n-1)$-dimensional Hausdorff measure $\mathcal H^{n-1}$ of $Z_{\phi_\lambda}$.

The main result is:
\begin{maintheo}\label{theorem} There exists a constant $c$ independent of $\lambda$ so that
\begin{equation}\label{LB}
\mathcal{H} ^{n-1} (Z_{\phi_\lambda}) \geq c\;
\left\{\begin{array}{ll}
\lambda^{\frac{3-n}{4}} \qquad \; n\leq 5 \\
\lambda^{\frac{17-5n}{16}} \qquad n >5
\end{array} \right. .
\end{equation}
\end{maintheo}
Hence in dimension $n=3$ we get a uniform lower bound for $ \mathcal{H} ^{n-1} (Z_{\phi_\lambda})$ which was proved recently by Colding-Minicozzi \cite{CM}. Here we give a new proof of their results for dimensions $n \leq 5$ but instead of local analysis of eigenfunctions we use the global methods of \cite{SZ}.

\subsection{Background}
The well known conjecture of S. T. Yau \cite{Y} states that there exist constants $c$ and $C$ independent of $\lambda$ such that
$$c \sqrt{\lambda} \leq \mathcal H ^{n-1} (Z_{\phi_\lambda}) \leq C \sqrt{\lambda}.$$
The conjecture was proved by Donnelly and Fefferman \cite{DF} when $(M,g)$ is real analytic.

In dimension $n=2$ and the $C^\infty$ case the best bounds are
$$c\sqrt{\lambda} \leq \mathcal H ^{1} (Z_{\phi_\lambda}) \leq C \lambda ^{3/4}.$$
The lower bound was proved by Br\"uning \cite{B} and Yau. See also \cite{Sa} where the constant $c$ is found explicitly in terms of $(M,g)$. The upper bound for $n=2$ was proved by Donnelly-Fefferman \cite{DF2} and also Dong \cite{D}.

For dimensions $n \geq 3$ the existing estimates are very far from the conjecture. In fact until recently the best bounds were the following:
$$c^{-\sqrt{\lambda}} \leq \mathcal H ^{n-1} (Z_{\phi_\lambda}) \leq \lambda^{C\sqrt{\lambda}}.$$
The lower bound was proved by Han-Lin \cite{HL} and the upper bound was proved by Hardt-Simon \cite{HS}. However, recently three papers by Sogge-Zelditch\cite{SZ}, 
Colding-Minicozzi \cite{CM} and Mangoubi \cite{M} were published where the lower bounds are improved from being exponentially decaying to being polynomially decaying as $\lambda \to \infty$:
\begin{equation}\label{SZ}\text{Sogge-Zelditch:} \qquad \quad \mathcal H ^{n-1} (Z_{\phi_\lambda}) \geq c \lambda^{(7-3n)/8}, \end{equation}
\begin{equation}\label{CM}\text{Colding-Mincozzi:} \qquad \quad \mathcal H ^{n-1} (Z_{\phi_\lambda})\geq c \lambda^{(3-n)/4},  \end{equation}
$$ \quad \text{Mangoubi:}\qquad \quad \qquad  \mathcal H ^{n-1} (Z_{\phi_\lambda})\geq c \lambda^{(3-n)/2 -\frac{1}{2n}} .$$
The result (\ref{CM}) is the best lower bound to this day which in particular gives a uniform lower bound in dimension $n=3$. However, the methods of these papers are different from each other. The paper \cite{SZ} uses global analysis of eigenfunctions such as the identity (\ref{DSZ1}) below and some $L^\infty$ and $L^1$ estimates for eigenfunctions, while \cite{CM} and \cite{M} are closer to \cite{DF} in spirit and use local analysis of eigenfunctions in balls of radius ${\lambda}^{-1/2}$. Our approach is closely related to \cite{SZ} and relies on an identity which was proved in that paper. It is interesting that both local and global methods give us the same lower bound in dimension $n=3$.

\section{Dong-Sogge-Zelditch formula}
We first note that the singular set of $\phi_\lambda$ which is defined by
$$\Sigma=Z_{\phi_\lambda} \cap \{\nabla \phi_\lambda =0\},$$ has finite $(n-2)$-dimensional Hausdorff measure and hence zero $(n-1)$-dimensional Hausdorff measure (see \cite{Ch} also \cite{H, HL}). Therefore $Z_{\phi_\lambda}$ admits a natural Riemannian hypersurface measure $dS_g$.

The Dong-Sogge-Zelditch formula (see \cite{D, SZ} and also \cite{ACF}) states that for every $f \in C^\infty(M)$:
\begin{equation}\label{DSZ} \int_M ((\Delta+\lambda)f)\, |\phi_\lambda|\, dV_g=2\int_{Z_{\phi_\lambda}} f\,|\nabla \phi_\lambda|\, dS_g.
\end{equation}
A special case of this formula was proved and used by Dong \cite{D} with $f= (|\nabla \phi_\lambda|^2 + \frac{\lambda^2\phi_\lambda^2}{n^2})^{-1/2}$ to obtain the upper bound $\lambda^{3/4}$ for the length of the nodal lines when $n=2$. The identity (\ref{DSZ}) was proved in \cite{SZ} and it was used with $f=1$ to obtain the lower bound (\ref{SZ}). To prove our theorem we will use their identity with $f=1$ and $f=|\nabla \phi_\lambda|^2$.

\subsection{Sogge-Zelditch method}  Since some estimates of \cite{SZ} will be used in our proof we explain their method in this section.

First they put $f=1$ in equation (\ref{DSZ}) to get
\begin{equation}\label{DSZ1}
\lambda \int_M |\phi_\lambda| dV_g=2\int_{Z_{\phi_\lambda}} |\nabla \phi_\lambda|\, dS_g.
\end{equation}
Then to find a lower bound for $\mathcal{H}^{n-1} (Z_{\phi_\lambda})=\int_{Z_{\phi_\lambda}} dS_g$ they find an upper bound for $||\nabla \phi_\lambda||_{L^\infty}$ and a lower bound for $|| \phi_\lambda||_{L^1}$:
\begin{equation}\label{grad}
||\nabla \phi_\lambda||_{L^\infty} \leq C\lambda^{(n+1)/4}
\end{equation} and
\begin{equation}\label{L1}
|| \phi_\lambda||_{L^1} \geq C \lambda^{(1-n)/8}.
\end{equation}
Applying these estimates to (\ref{DSZ1}) gives (\ref{SZ}). Both of these estimates are sharp. The upper bound (\ref{grad}) for $||\nabla \phi_\lambda||_{L^\infty}$ is achieved by the zonal spherical harmonics on $S^2$ and the highest weight spherical harmonics (Gaussian beams) on $S^2$ saturate the lower bound (\ref{L1}) for $|| \phi_\lambda||_{L^1}$.
The estimate (\ref{grad}) is proved using a local Weyl law. One way to prove (\ref{L1}) is to do the following:

By the H\"older inequality
$$ 1=||\phi_\lambda||_{L^2}\leq  ||\phi_\lambda||^{\frac{p-2}{2(p-1)}}_{L^1} ||\phi_\lambda||^{\frac{p}{2(p-1)}}_{L^p},$$
therefore
$$|| \phi_\lambda||_{L^1} \geq C || \phi_\lambda||^{-\frac{p}{p-2}}_{L^p}.$$
Then (\ref{L1}) follows from Sogge $L^p$-estimates \cite{So}:
\begin{equation}\label{Sogge}
||\phi_\lambda||_{L^p} \leq C \lambda^{\delta(p)} \quad \text{where} \quad  \delta(p)= \left\{\begin{array}{ll} \frac{n(p-2)-p}{4p} \qquad p\geq \frac{2(n+1)}{n-1} \\
 \frac{(n-1)(p-2)}{8p} \quad 2\leq p \leq \frac{2(n+1)}{n-1}  \end{array} \right. ,
\end{equation}
with $p= \frac{2(n+1)}{n-1}$.

\section{Proof of the theorem}
To prove our theorem we first apply the H\"older inequality to (\ref{DSZ1}):
\begin{equation}\label{holder}
\lambda \int_M |\phi_\lambda| dV_g=2\int_{Z_{\phi_\lambda}} |\nabla \phi_\lambda|\, dS_g \leq 2\, (\mathcal{H}^{n-1} (Z_{\phi_\lambda}))^{\frac{2}{3}}
\; \big(\int_{Z_{\phi_\lambda}} | \nabla \phi_\lambda|^3 \, dS_g \big)^{\frac{1}{3}} .
\end{equation}
To estimate $|| \nabla \phi_\lambda||_{L^3(Z_{\phi_\lambda})}$ we use (\ref{DSZ}) with $f=|\nabla \phi_\lambda|^2$ and we find that
\begin{equation}\label{L3}
\int_{Z_{\phi_\lambda}} | \nabla \phi_\lambda|^3 dS_g=\frac{1}{2} \int_M \big((\Delta+\lambda)|\nabla \phi_\lambda|^2\big) |\phi_\lambda|\,dV_g.
\end{equation}
The Bochner identity for $\frac{1}{2}\Delta |\nabla \phi_\lambda|^2$ is
\begin{equation}\label{Bochner}
\frac{1}{2} \Delta |\nabla \phi_\lambda|^2= \;|H(\phi_\lambda)|^2-\lambda |\nabla \phi_\lambda|^2+
\; \text{Ric}\,(\nabla \phi_\lambda, \nabla \phi_\lambda),
\end{equation}
where $H(\phi_\lambda)$ is the Hessian of $\phi_\lambda$ and $|H(\phi_\lambda)|=\,\sqrt{\text{Tr}\,(H(\phi_\lambda)^2)}$ is the standard Riemannian norm of $H(\phi_\lambda)$.

By applying (\ref{Bochner}) to (\ref{L3}) we obtain
\begin{equation}\label{L3bochner}
\int_{Z_{\phi_\lambda}}| \nabla \phi_\lambda|^3dS_g=\int_M \big\{ \; |H(\phi_\lambda)|^2 |\phi_\lambda|-\frac{1}{2}\lambda |\nabla \phi_\lambda|^2 |\phi_\lambda| + \;
 \text{Ric}\,(\nabla \phi_\lambda, \nabla \phi_\lambda) |\phi_\lambda|\big\} dV_g.
\end{equation}
Since $M$ is compact we can find a constant $C$ such that $$  \text{Ric}\,(\nabla \phi_\lambda, \nabla \phi_\lambda)
\leq C |\nabla \phi_\lambda|^2.$$ Thus in (\ref{L3bochner}) for $\lambda$ large enough the sum of the second and third terms is negative and can be neglected in our estimates. We then bound the first term in (\ref{L3bochner}) as follows:
\begin{equation}\label{bound}\int_M  |H(\phi_\lambda)|^2 |\phi_\lambda| dV_g\leq || H(\phi_\lambda)||^{2}_{L^3}||\phi_\lambda||_{L^3} \leq C \lambda^2 || \phi_\lambda ||^3_{L^3} \leq C
\left\{\begin{array}{ll} \lambda^{2+\frac{n-1}{8}} \qquad n\leq 5\\
\lambda^{2+\frac{n-3}{4}} \quad \quad n >5  \end{array} \right., \end{equation}
where we have used the H\"older inequality in the first, elliptic $W^{2,p}$ estimates in the second inequality and Sogge $L^3$-estimates (\ref{Sogge}) in the last inequality. The theorem follows if we combine (\ref{bound}), (\ref{L3bochner}), (\ref{holder}) and (\ref{L1}).

We close this section with two remarks on our proof.

\begin{rem} In the proof of
$$|| H(\phi_\lambda)||^2_{L^3} \leq C \lambda^2 || \phi_\lambda ||^2_{L^3} $$ we have applied the following elliptic estimates (see \cite{GT}, Theorem $9.11$) to the equation $\Delta \phi_\lambda = -\lambda \phi_\lambda$ (when it is written in a local chart):

Let $V \subset \mathbb R^n$ and open set and Let $\phi \in C^\infty(V)$ satisfy $L \phi =g$ where $L$ is a second order elliptic differential operator with smooth coefficients. Then for any domain $U \subset \subset V$:
$$ || \phi|| _{W^{2,p}(U)} \leq C (||g ||_{L^p(V)}+|| \phi||_{L^p(V)}), \qquad 1 < p < \infty,$$
where $C$ depends only on $n$, $p$, $U$, $V$ and the moduli of continuity of the coefficients of $L$.
\end{rem}

\begin{rem}We may lose too much when we ignore the negative term $-\frac{1}{2} \lambda\int_M |\nabla \phi_\lambda|^2 |\phi_\lambda|$  in (\ref{L3bochner}). Any possible cancellation between the first and the second terms in (\ref{L3bochner}) which reduces the order of $\lambda$ will improve the existing results on the lower bounds of the volumes of nodal sets.
\end{rem}

\begin{rem}In (\ref{holder}) we can use $||\nabla \phi_\lambda||_{L^p(Z_{\phi_\lambda})}$ for different powers of $p$ and use (\ref{DSZ}) with $f=|\nabla \phi_\lambda|^{p-1}$. However it seems if we follow this method for $p\geq 3$ we find that the power $p$ which works the best is $p=3$. Of course we can use $p<3$ but this will make the integrand on the left hand side of (\ref{DSZ}) a singular function which we do not know how to treat.
\end{rem}

\section{Methods of Donnelly-Fefferman,  Colding-Minicozzi and Han-Lin}
In this section we explain the main idea of the paper \cite{CM} and the closely related section $6.2$ of the book \cite{HL} and at the end we will list a couple of related questions.

 The main idea which goes back to \cite{DF} is to cover the manifold with balls $B$ of radius $\frac{a}{\sqrt{\lambda}}$ where $a$ is chosen large enough so that each ball contains a zero of $\phi_\lambda$. See Lemma $6.2.1$ of \cite{HL} for the existence of such $a$. By choosing a larger $a$ we can assume that the balls are centered at a zero of $\phi_\lambda$. Then a ball $B$ (which is always assumed to have radius $\frac{a}{\sqrt{\lambda}}$) is called good if the following doubling estimate holds for a constant $K$ independent of $\lambda$:
\begin{equation}\label{doubling}
\int_{2B} |\phi_\lambda|^2 dV_g \leq K \int_B |\phi_\lambda|^2 dV_g
\end{equation}
We note that for a general ball the best doubling estimate is far from being good and states that $K \leq c^{\-\sqrt{\lambda}}$ (see for example Lemma $6.1.1$ of \cite{HL}). If $B$ is a good ball then one can prove
\begin{equation}\label{localvol}\mathcal{H}^{n-1} (Z_{\phi_\lambda} \cap B)\geq C \lambda^{-(n-1)/2} .\end{equation}
For a proof of this see Lemma $6.2.4$ of \cite{HL} (see also Proposition $1$ of \cite{CM}). So the main question to ask in finding a lower bound for
$\mathcal{H}^{n-1} (Z_{\phi_\lambda})$ is ``how many good balls are there?'' In \cite{CM} the authors use Sogge $L^p$-estimates and show that there are at least $ \lambda^{(n+1)/4}$ good balls hence they obtain the lower bound $\lambda^{(n+1)/4} \lambda ^{-(n-1)/2}=\lambda^{(3-n)/4}$.

Now we list a couple of remarks and questions regarding this local method:
\begin{itemize}
\item[1.] The above approach was motivated by the work of Donnelly-Fefferman \cite{DF} where a good ball was defined similarly except that instead of $L^2$ doubling estimates $L^\infty$ doubling estimates were considered. They showed that in the real analytic case the total volume of good balls is proportional to the volume of $(M,g)$ hence the total number of good balls is of size $\lambda^{n/2}$.  Using  (\ref{localvol}) this implies that $\mathcal{H}^{n-1} (Z_{\phi_\lambda}) \geq c \sqrt{\lambda}$.
    \\

\item[2.] For $n=2$ the proof of Br\"uning for the lower bound $c\sqrt{\lambda}$ is similar to the local argument above but no doubling estimates (or notion of good balls) are needed. First we choose $B$ of radius $\frac{a}{\sqrt{\lambda}}$ as above such that it is centered at a zero of $\phi_\lambda$. Then we shrink $B$ to a smaller ball (still called $B$) of radius  $\frac{\epsilon}{\sqrt{\lambda}}$ and with the same center ($\epsilon$ to be chosen later). These balls do no necessarily cover $M$ anymore but we have not changed the number of balls which is roughly $\lambda^{n/2}$.  In dimension $n=2$ one can prove directly that (\ref{localvol}) is true for every ball $B$ of radius $\frac{\epsilon}{\sqrt{\lambda}}$. To do this we rescale the ball $B$ to the Euclidean ball $\widetilde{B}$ of radius one so the operator $\Delta+\lambda$ becomes a small perturbation of the Laplacian. The constant $\epsilon$ is chosen small enough so that the maximum principle applies. Then by the maximum principle the rescaled eigenfunction does not have any loops (closed curves) in its set of nodal lines in $\widetilde{B}$. Since it has a zero at the origin therefore there is a nodal curve of constant length in $\widetilde{B}$ and therefore a nodal curve of length $c\lambda^{-1/2}$ in $B$.
\\

\item[3.] It is clear from (\ref{localvol}) that if we have enough number of good balls then we get the lower bound in Yau's conjecture. One hopes to prove (\ref{doubling}) for a quantum ergodic sequence of eigenfunctions where eigenfunctions distribute on open sets according to the volume of open sets. But a ball of radius $\lambda^{-1/2}$ seems too small for semiclassical techniques to be applied.
\\

\item[4.] If we follow the paper of \cite{CM} it is clear that if we have an orthonormal basis (or any sequence) of eigenfunctions which is uniformly bounded in $L^\infty$ then we get the lower bound $\mathcal{H}^{n-1} (Z_{\phi_\lambda}) \geq c \sqrt{\lambda}$. So far the only example of a Riemannian manifold for which we know there exists a uniformly bounded orthonormal basis of eigenfunctions is the flat torus. For example we do not even know if such a basis exists for the round sphere $S^2$. A random wave (as a model for quantum ergodic eigenfunctions) on $S^2$ has the $L^\infty$ bound $\sqrt{\log \lambda}$.
\end{itemize}

\subsection{Acknowledgments:} The authors are grateful to Christopher Sogge and Steve Zelditch for their generous paper \cite{SZ} which not only includes detailed proofs of their results but also provides many examples and methods for possible improvements. The first author would like to thank Steve Zelditch, Maciej Zworski and Colin Guillarmou for helpful discussions and Kiril Datchev for reading the earlier version of this paper.

\end{document}